\documentclass[a4paper,12pt,eqno]{article}
\usepackage{theorem}
\usepackage{latexsym,amssymb,amsfonts,amsmath}
\usepackage{graphicx}
\usepackage{amsmath}

\newtheorem{thm}{Theorem}[section]

\newtheorem{cor}[thm]{Corollary}
\newtheorem{pro}[thm]{Proposition}

\theoremheaderfont{\scshape}
\newcommand{\RM}{\mathbb{R}}

\title{{\Large {\bf VERTEX-FACE/ZETA CORRESPONDENCE}}
\author{
{\small Takashi Komatsu} \\
{\scriptsize Math. Research Institute Calc for Industry} \\
{\scriptsize Minami, Hiroshima, 732-0816, Japan} \\ 
{\scriptsize e-mail: ta.komatsu@sunmath-calc.co.jp} \\ 
{\small Norio Konno}\\
{\scriptsize Department of Applied Mathematics, 
Faculty of Engineering, 
Yokohama National University}\\
{\scriptsize Hodogaya, Yokohama, 240-8501, Japan}\\
{\scriptsize e-mail: konno-norio-bt@ynu.ac.jp}\\
{\small Iwao Sato}\\
{\scriptsize Oyama National College of Technology}\\
{\scriptsize Oyama, Tochigi 323-0806, Japan}\\
{\scriptsize e-mail: isato@oyama-ct.ac.jp}\\}
}
\date{\empty }
\begin{document}
\maketitle

\par\noindent
\begin{small}
\par\noindent
{\bf Abstract}. 
We present the characteristic polynomial for the transition matrix of a vertex-face walk on a graph, and obtain its spectra. 
Furthermore, we express the characteristic polynomial for the transition matrix of a vertex-face walk on the 2-dimensional torus 
by using its adjacency matrix, and obtain its spectra. 
As an application, we define a new walk-type zeta function with respect to the transition matrix of a vertex-face walk 
on the 2-dimensional torus, and present its explicit formula.

\footnote[0]{
{\it Abbr. title:} Vertex-face/zeta correspondence
}
\end{small}

\section{Introduction} 

Recently, there were exciting developments between quantum walk \cite{Ambainis2003, Kempe2003, 
Kendon2007, Konno2008b, VA} on a graph and the Ihara zeta function \cite{Ihara1966, Serre, Sunada1986, Sunada1988, Hashimoto1989, Bass1992, 
ST1996, FZ1999, KS2000} 
of a graph: 
Grover/Zeta correspondence \cite{KomatsuEtAl2021}; Walk/Zeta correspondence \cite{KomatsuEtAl2021b}. 

In Grover/Zeta correspondence \cite{KomatsuEtAl2021}, a zeta function and a generalized zeta function of a graph $G$ 
with respect to its Grover matrix as analogue of the Ihara zeta function and 
the generalized Ihara zeta function \cite{ChintaEtAl} of $G$ were defined. 
By using the Konno-Sato theorem \cite{KS2011}, the limits on the generalized zeta functions and the generalized Ihara zeta functions 
of a family of finite regular graphs were written as an integral expression, and contained the result on the generalized Ihara zeta function 
in Chinta et al. \cite {ChintaEtAl}. 
Furthermore, the limit on the generalized Ihara zeta functions of a family of finite torus is written as an integral expression, 
and contained the result on the Ihara zeta function of the two-dimensional integer lattice $\mathbb{Z}^2$ in Clair \cite {Clair}. 

In Walk/Zeta correspondence \cite{KomatsuEtAl2021b}, a walk-type zeta function was defined without use of the determinant expressions of zeta function 
of a graph $G$, and various properties of walk-type zeta functions of random walk (RW), correlated random walk (CRW) and quantum walk (QW) on $G$. 
Also, their limit formulas by using integral expressions were presented. 

Recently, Zhang \cite{Zhan} introduced a vertex-face walk on an orientable embedding of a graph, and presented the spectra for its transition matrix. 
 
In this paper, we treat a walk-type zeta function of a vertex-face walk on a graph defined by Zhang \cite{Zhan}. 

The rest of the paper is organized as follows. 
Section 2 gives a short review for Grover/Zeta correspondence. 
In Sect. 3, we state Walk/Zeta correspondence on finite torus. 
In Sect. 4, we present an explicit formula for the characteristic polynomial of the transition matrix of a vertex-face walk on a graph, 
and obtain its spectra. 
In Sect. 5, we express the characteristic polynomial for the transition matrix of a vertex-face walk on the 2-dimensional torus 
by using its adjacency matrix, and obtain its spectra. 
As an application, we define a new walk-type zeta function with respect to the transition matrix of a vertex-face walk on the 2-dimensional torus, 
and present its explicit formula.

\section{Grover/Zeta correspondence}

All graphs in this paper are assumed to be simple. 
Let $G=(V(G),E(G))$ be a connected graph (without multiple edges and loops) 
with the set $V(G)$ of vertices and the set $E(G)$ of unoriented edges $uv$ 
joining two vertices $u$ and $v$.
Furthermore, let $n=|V(G)|$ and $m=|E(G)|$ be the number of vertices and edges of $G$, respectively. 
For $uv \in E(G)$, an arc $(u,v)$ is the oriented edge from $u$ to $v$. 
Let $D_G$ the symmetric digraph corresponding to $G$. 
Set $D(G)= \{ (u,v),(v,u) \mid uv \in E(G) \} $. 
For $e=(u,v) \in D(G)$, set $u=o(e)$ and $v=t(e)$. 
Furthermore, let $e^{-1}=(v,u)$ be the {\em inverse} of $e=(u,v)$. 
For $v \in V(G)$, the {\em degree} $\deg {}_G \ v = \deg v = d_v $ of $v$ is the number of vertices 
adjacent to $v$ in $G$.  

A {\em path $P$ of length $n$} in $G$ is a sequence 
$P=(e_1, \ldots ,e_n )$ of $n$ arcs such that $e_i \in D(G)$,
$t( e_i )=o( e_{i+1} )(1 \leq i \leq n-1)$. 
If $e_i =( v_{i-1} , v_i )$ for $i=1, \cdots , n$, then we write 
$P=(v_0, v_1, \ldots ,v_{n-1}, v_n )$. 
Set $ \mid P \mid =n$, $o(P)=o( e_1 )$ and $t(P)=t( e_n )$. 
Also, $P$ is called an {\em $(o(P),t(P))$-path}. 
We say that a path $P=( e_1 , \ldots , e_n )$ has a {\em backtracking} 
if $ e^{-1}_{i+1} =e_i $ for some $i(1 \leq i \leq n-1)$. 
A $(v, w)$-path is called a {\em $v$-cycle} 
(or {\em $v$-closed path}) if $v=w$. 
Let $B^r$ be the cycle obtained by going $r$ times around a cycle $B$. 
Such a cycle is called a {\em multiple} of $B$. 
A cycle $C$ is {\em reduced} if 
both $C$ and $C^2 $ have no backtracking. 

The {\em Ihara zeta function} of a graph $G$ is 
a function of a complex variable $u$ with $|u|$ 
sufficiently small, defined by 
\[
{\bf Z} (G, u)= \exp \left( \sum^{\infty}_{k=1} \frac{N_k }{k} u^k \right) , 
\]
where $N_k $ is the number of reduced cycles of length $k$ in $G$.

Let $G$ be a connected graph with $n$ vertices $v_1, \ldots ,v_n $. 
The {\em adjacency matrix} ${\bf A}= {\bf A} (G)=(a_{ij} )$ is 
the square matrix such that $a_{ij} =1$ if $v_i$ and $v_j$ are adjacent, 
and $a_{ij} =0$ otherwise.
If $ \deg {}_G \ v=k$(constant) for each $v \in V(G)$, then $G$ is called 
{\em $k$-regular}.

\begin{thm}[Ihara; Bass] 
Let $G$ be a connected graph. 
Then the reciprocal of the Ihara zeta function of $G$ is given by 
\[
{\bf Z} (G,u )^{-1} =(1- u^2 )^{r-1} 
\det ( {\bf I} -u {\bf A} (G)+ u^2 ( {\bf D} - {\bf I} )) , 
\]
where $r$ is the Betti number of $G$, 
and ${\bf D} =( d_{ij} )$ is the diagonal matrix 
with $d_{ii} = \deg v_i$ and $d_{ij} =0, i \neq j , 
(V(G)= \{ v_1 , \ldots , v_n \} )$. 
\end{thm}

Let $G=(V(G),E(G))$ be a connected graph with $\nu $ vertices and $ x_0 \in V(G)$ 
a fixed vertex.  
Then the {\em generalized Ihara zeta function} $\zeta {}_G (u)$ of $G$ is defined by 
\[
\zeta {}_G (u)= \zeta (G, u)= \exp \left( \sum^{\infty}_{m=1} \frac{N^0_m }{m} u^m \right) , 
\]
where $N^0_m $ is the number of reduced $x_0$-cycles of length $m$ in $G$. 
Furthermore, the {\em Laplacian} of $G$ is given by 
\[
\Delta {}_{\nu } = \Delta (G) = {\bf D} - {\bf A} (G). 
\]

A formula for the generalized Ihara zeta function of a vertex transitive graph is given as follows:

\begin{thm}[Chinta, Jorgenson and Karlsson] 
Let $G$ be a  vertex-transitive $(q+1)$-regular graph with spectral measure $\mu {}_{\Delta }$ for the Laplacian $\Delta $. 
Then 
\[
\zeta {}_G (u)^{-1} =(1-u^2 )^{(q-1)/2} \exp ( \int \log (1-(q+1- \lambda )u+q u^2 ) d \mu {}_{\Delta } ( \lambda )) . 
\]
\end{thm}

A graph $G$ is called {\em vertex-transitive} if there exists an automorphism $ \phi $ of the automorphism group 
$Aut \ G$ of $G$ such that $ \phi (u)=v$ for each $u,v \in V(G)$. 
Note, if $G$ is a vertex-transitive graph with $n$ vertices, then 
\[
\zeta {}_G (u)= {\bf Z} (G,u)^{1/n} . 
\]

Let $G$ be a connected graph with $n$ vertices and $m$ edges. 
Set $V(G)= \{ v_1 , \ldots , v_n \} $ and $d_j = d_{v_j} = \deg v_j , \ j=1, \ldots , n$. 
Then the {\em Grover matrix} ${\bf U} ={\bf U} (G)=( U_{ef} )_{e,f \in D(G)} $ 
of $G$ is defined by 
\[
U_{ef} =\left\{
\begin{array}{ll}
2/d_{t(f)} (=2/d_{o(e)} ) & \mbox{if $t(f)=o(e)$ and $f \neq e^{-1} $, } \\
2/d_{t(f)} -1 & \mbox{if $f= e^{-1} $, } \\
0 & \mbox{otherwise. }
\end{array}
\right. 
\]
The discrete-time quantum walk with the matrix ${\bf U} $ as a time evolution matrix 
is called the {\em Grover walk} on $G$.

Let $G$ be a connected graph with $\nu$ vertices and $m$ edges. 
Then the $\nu \times \nu$ matrix ${\bf P}_{\nu } = {\bf P} (G)=( P_{uv} )_{u,v \in V(G)}$ is given as follows: 
\[
P_{uv} =\left\{
\begin{array}{ll}
1/( \deg {}_G \ u)  & \mbox{if $(u,v) \in D(G)$, } \\
0 & \mbox{otherwise.}
\end{array}
\right.
\] 
Note that the matrix ${\bf P} (G)$ is the transition probability matrix of the simple random walk on $G$.

We introduce the {\em positive support} ${\bf F}^+ =( F^+_{ij} )$ of 
a real matrix ${\bf F} =( F_{ij} )$ as follows: 
\[
F^+_{ij} =\left\{
\begin{array}{ll}
1 & \mbox{if $F_{ij} >0$, } \\
0 & \mbox{otherwise. }
\end{array}
\right.
\] 
Ren et al. \cite{RenETAL} showed that the edge matrix of a graph is the positive support $({\bf U}^T )^+ $ 
of the transpose of its Grover matrix ${\bf U} $, i.e., 
\[
{\bf Z} (G,u)^{-1} = \det ( {\bf I}_{2m} -u {\bf U}^+ ) . 
\]
The Ihara zeta function of a graph is just a zeta function on the positive support of 
the Grover matrix of a graph.

Note that, by Theorem 2.1,  
\[
\det ( {\bf I}_{2m} -u {\bf U}^+ )=(1-u^2)^{m- \nu} \det ((1+qu^2) {\bf I}_{\nu} -((q+1) {\bf I}_{\nu } - \Delta_{\nu } )u) . 
\]

Now, we propose a new zeta function of a graph.  
Let $G$ be a connected graph with $m$ edges. 
Then we define a zeta function $ \overline{{\bf Z}} (G, u)$ of $G$ satisfying   
\[
\overline{{\bf Z}} (u)^{-1} =\overline{{\bf Z}} (G, u)^{-1} = \det ( {\bf I}_{2m} -u {\bf U} ) .    
\]

In Konno and Sato \cite{KS2011}, they presented the following results in our setting.

\begin{thm}[Konno and Sato]
Let $G$ be a connected graph with $\nu$ and $m$ edges. 
Then  
\[
\det ( {\bf I}_{2m} -u {\bf U} )=(1-u^2)^{m-\nu} \det ((1+u^2) {\bf I}_{\nu} -2u {\bf P} (G)) . 
\]
\end{thm}

We give a weight functions $w: D(G) \times D(G) \longrightarrow \mathbb{C} $ as follows: 
\[
w (e,f) =\left\{
\begin{array}{ll}
2/ \deg t(e)  & \mbox{if $t(e)=o(f)$ and $f \neq e^{-1} $, } \\
2/ \deg t(e) -1 & \mbox{if $f= e^{-1} $, } \\
0 & \mbox{otherwise. }
\end{array}
\right.
\] 
For a cycle $C=( e_1, e_2 , \ldots , e_r )$, let 
\[
w(C)=w(e_1 , e_2 ) \cdots w( e_{r-1} , e_r) w( e_r , e_1 ) . 
\]

We define a generalized zeta function with respect to the Grover matrix of a graph. 
Let $G=(V(G),E(G))$ be a connected graph and $ x_0 \in V(G)$ a fixed vertex.  
Then the {\em generalized zeta function} $\overline{\zeta} {}_G (u)$ of $G$ is defined by 
\[
\overline{\zeta} {}_G (u)= \overline{\zeta} (G, u)=\exp \left( \sum^{\infty}_{r=1} \frac{N^0_r }{r} u^r \right) , 
\]
where 
\[
N^0_r = \sum \{ w(C) \mid C: \ an \ x_0-cycle \ of \ length \ r \ in \ G \} . 
\] 
Note, if $G$ is a vertex-transitive graph with $n$ vertices, then 
\begin{equation} 
\overline{\zeta} {}_G (u)= \overline{{\bf Z}} (G,u)^{1/n} . 
\end{equation}

Then we obtain the following results for a series of finite vertex-transitive $(q+1)$-regular graphs.

\begin{thm}[Grover/Zeta correspondence]  
Let $\{ G_n \}^{\infty}_{n=1} $ be a series of finite vertex-transitive $(q+1)$-regular graphs such that 
\[
\lim {}_{n \rightarrow \infty} |V(G_n )|= \infty . 
\]
Then
\begin{enumerate}   
\item $\lim_{n \rightarrow \infty} \overline{\zeta}_{G_n} (u)^{-1} =
(1-u^2 )^{(q-1)/2} \exp [ \int \log \{ ((1+u^2 )-2u \lambda )) \} d \mu_{P} ( \lambda )] $,   
\item $\lim_{n \rightarrow \infty} \overline{\zeta}_{G_n } (u)^{-1} =
(1-u^2 )^{(q-1)/2} \exp [ \int \log \{ (1-2u+u^2 )+ \frac{2u}{q+1} \lambda \} d \mu_{\Delta} ( \lambda )] $, 
\item $\lim_{n \rightarrow \infty} \zeta_{G_n} (u)^{-1} =
(1-u^2 )^{(q-1)/2} \exp \left[ \int \log \{ (1+q u^2 )-(q+1)u \lambda \} d \mu_P ( \lambda ) \right] $,    
\item $\lim_{n \rightarrow \infty} \zeta_{G_n} (u)^{-1} =
(1-u^2 )^{(q-1)/2} \exp \left[ \int \log \{ (1+q u^2 )-((q+1)- \lambda ) u \} d \mu_{\Delta } ( \lambda ) \right] $,  
\end{enumerate} 
where $d \mu_P ( \lambda )$ and $d \mu_{\Delta } ( \lambda )$ are the spectral measures for the transition operator ${\bf P}$ 
and the Laplacian $\Delta $. 
\end{thm}

We should note that the fourth formula in Theorem 2.4 is nothing but Theorem 1.3 in Chinta et al. [3].

Next, we obtain the following results for the generalized zeta function and 
the generalized Ihara zeta function of the $d$-dimensional integer lattice $\mathbb{Z}^d \ (d \geq 2)$.  

Let $T^d_N \ (d \geq 2)$ be the {\em $d$-dimensional torus (graph)} with $N^d$ veritices. 
Its vertices are located in coordinates $i_1 , i_2 , \ldots , i_d $ of a $d$-dimensional Euclidian space $\mathbb{R}^d $, 
where $i_j \in \{ 0,1, \ldots , N-1 \} $ for any $j$ from 1 to $d$. 
A vertex $v$ is adjacent to a vertex $w$ if and only if they have $d-1$ coordinates that are the same, 
and for the remaining coordinate $k$, we have $|i^v_k - i^w_k |=1$, where $i^v_k $ and $i^w_k $ are the $k$-th coordinate 
of $v$ and $w$, respectively. 
Then we have 
\[
|E( T^d_N )|=d N^d ,  
\]
and $T^d_N$ is a vertex-transitive $2d$-regular graph.

\begin{thm}[Grover/Zeta correspondence($T^d_N$ case )]    
Let $T^d_N \ (d \geq 2)$ be the $d$-dimensional torus with $N^d$ veritices. 
Then  
\[ 
\displaystyle 
\lim_{n \rightarrow \infty} \overline{\zeta} (T^d_N ,u)^{-1} =(1- u^2 )^{d-1} \exp \left[ \int^{2 \pi}_{0} \dots \int^{2 \pi}_{0}  
\log \{ (1+u^2 )- \frac{2u}{d} \sum^d_{j=1} \cos \theta_j \} \frac{d \theta_1}{2 \pi } \cdots \frac{d \theta_d}{2 \pi } \right] ,  
\]
\[ 
\displaystyle 
\lim_{n \rightarrow \infty} {\zeta} (T^d_N ,u)^{-1} =(1- u^2 )^{d-1} \exp \left[ \int^{2 \pi}_{0} \dots \int^{2 \pi}_{0}  
\log \{ (1+(2d-1) u^2 )-2u \sum^d_{j=1} \cos \theta_j \} \frac{d \theta_1}{2 \pi } \cdots \frac{d \theta_d}{2 \pi } \right] ,  
\] 
where $\int^{2 \pi}_{0} \dots \int^{2 \pi}_{0} $ is the $d$-th multiple integral and 
$ \frac{d \theta_1}{2 \pi } \cdots \frac{d \theta_d}{2 \pi } $ is the uniform measure on $[0, 2 \pi )^d $.  
\end{thm}

Specially, in the case of $d=2$, we obtain the following result.

\begin{cor}     
Let $T^2_N $ be the $2$-dimensional torus with $N^2$ veritices. 
Then  
\[  
\displaystyle 
\lim_{n \rightarrow \infty} \overline{\zeta} (T^2_N ,u)^{-1} =(1- u^2 ) \exp \left[ \int^{2 \pi}_{0} \int^{2 \pi}_{0}  
\log \{ (1+u^2 )-u \sum^d_{j=1} \cos \theta_j \} \frac{d \theta_1}{2 \pi } \frac{d \theta_2}{2 \pi } \right] ,  
\] 
\[ 
\displaystyle 
\lim_{n \rightarrow \infty} {\zeta} (T^2_N ,u)^{-1} =(1- u^2 ) \exp \left[ \int^{2 \pi}_{0} \int^{2 \pi}_{0}  
\log \{ (1+3 u^2 )-2u \sum^2_{j=1} \cos \theta_j \} \frac{d \theta_1}{2 \pi } \frac{d \theta_2}{2 \pi } \right] .   
\] 
\end{cor}

The second formula corresponds to Equation (10) in Clair [4].

\section{Walk/Zeta correspondence on torus} 

We state the Walk/Zeta correspondence on a finite torus. 
At first, we give the definition of the $2d$-state discrete-time walk on $T^d_N$. 
The discrete-time walk is defined by using a {\em shift operator} and a {\em coin matrix} which will be mentioned below.

Let $f : T^d_N \longrightarrow \mathbb{C}^{2d}$.
For $j = 1,2,\ldots,d$ and ${\bf x} \in T^d_N$, the shift operator $\tau_j$ is defined by 
\begin{align*}
(\tau_j f)({\bf x}) = f( {\bf x} - {\bf e}_{j}),
\end{align*} 
where $\{ {\bf e}_1,{\bf e}_2,\ldots,{\bf e}_d \}$ denotes the standard basis of $\mathbb{R}^d$.

Let $A=[a_{ij}]_{i,j=1,2,\ldots,2d}$ be a $2d \times 2d$ matrix with $a_{ij} \in \mathbb{C}$ for $i,j =1,2,\ldots,2d$. 
We call $A$ the {\em coin matrix}. If $a_{ij} \in [0,1]$ and $\sum_{i=1}^{2d} a_{ij} = 1$ for any $j=1,2, \ldots, 2d$, 
then the walk is a CRW. 
In particular, when $a_{i1} = a_{i2} = \cdots = a_{i 2d}$ for any $i=1,2, \ldots, 2d$, this CRW becomes a RW. If $A$ is unitary, 
then the walk is a QW. So our class of walks contains RWs, CRWs, and QWs as special models.

To describe the evolution of the walk, we decompose the $2d \times 2d$ coin matrix $A$ as
\begin{align*}
A=\sum_{j=1}^{2d} P_{j} A,
\end{align*}
where $P_j$ denotes the orthogonal projection onto the one-dimensional subspace $\mathbb{C}\eta_j$ in $\mathbb{C}^{2d}$. 
Here $\{\eta_1,\eta_2, \ldots, \eta_{2d}\}$ denotes a standard basis on $\mathbb{C}^{2d}$.
 
The discrete-time walk associated with the coin matrix $A$ on $T^d_N$ is determined by the $2d N^d \times 2d N^d$ matrix
\begin{align}
M_A=\sum_{j=1}^d \Big( P_{2j-1} A \tau_{j}^{-1} + P_{2j} A \tau_{j} \Big).
\label{unitaryop1}
\end{align}
Let $ \mathbb{Z}_{\geq} = \mathbb{Z} \cup \{ 0 \} $. 
Then the state at time $n \in \mathbb{Z}_{\ge}$ and location ${\bf x} \in T^d_N$ can be expressed by a $2d$-dimensional vector:
\begin{align*}
\Psi_{n}( {\bf x} )=
\begin{bmatrix}
\Psi^{1}_{n}( {\bf x} ) \\ \Psi^{2}_{n}( {\bf x} ) \\ \vdots \\ \Psi^{2d}_{n}( {\bf x} ) 
\end{bmatrix} 
\in \mathbb{C}^{2d}.
\end{align*} 
For $\Psi_n : T^d_N \longrightarrow \mathbb{C}^{2d} \ (n \in \mathbb{Z}_{\geq})$, from Eq. (2), the evolution of the walk is defined by 
\begin{align}
\Psi_{n+1}( {\bf x} ) \equiv (M_A \Psi_{n})( {\bf x} )=\sum_{j=1}^{d}\Big(P_{2j-1}A\Psi_{n}( {\bf x} + {\bf e}_j)+P_{2j}A\Psi_{n}( {\bf x} - {\bf e}_j)\Big).
\end{align}

Now, we define the {\em walk-type zeta function} by 
\begin{align}
\overline{\zeta} \left(A, T^d_N, u \right) = \det \Big( I_{2d N^d} - u M_A \Big)^{-1/N^d}.
\label{satosan01}
\end{align}
In general, for a $d_c \times d_c$ coin matrix $A$, we put  
\begin{align*}
\overline{\zeta} \left(A, T^d_N, u \right) = \det \Big(I_{d_c N^d} - u M_A \Big)^{-1/N^d}.
\end{align*}

Komatsu et al. \cite{KomatsuEtAl2021b} obtained the following result. 
\begin{thm}
\begin{align*}
\overline{\zeta} \left(A, T^d_N, u \right) ^{-1}
&= \exp \left[ \frac{1}{N^d} \sum_{\widetilde{ {\bf k} } \in \widetilde{\mathbb{K}}_N^d} \log \left\{ \det \Big( F(\widetilde{ {\bf k} }, u) \Big) \right\} \right],
\\
\lim_{N \to \infty} \overline{\zeta} \left(A, T^d_N, u \right) ^{-1}
&=
\exp \left[ \int_{[0,2 \pi)^d} \log \left\{ \det \Big( F \left( \Theta^{(d)}, u \right)  \Big) \right\} d \Theta^{(d)}_{unif} \right],
\end{align*}
where $ \tilde{\mathbb{K}}_N = \{ 0, 2 \pi /N, \ldots , 2 \pi (N-1)/N \} $, $\Theta^{(d)} = (\theta_1, \theta_2, \ldots, \theta_d) (\in [0, 2 \pi)^d)$ 
and $d \Theta^{(d)}_{unif}$ denotes the uniform measure on $[0, 2 \pi)^d$, that is,
\begin{align*}
d \Theta^{(d)}_{unif} = \frac{d \theta_1}{2 \pi } \cdots \frac{d \theta_d}{2 \pi } .  
\end{align*}
Furthermore,  
\begin{align*}
F \left( {\bf w} , u \right) = I_{2d} - u \widehat{M}_A ( {\bf w} ) \ \ \  and \ \ \                           
\widehat{M}_A( {\bf w} )=\sum_{j=1}^{d} \Big( e^{i w_j } P_{2j-1} A + e^{-i w_j } P_{2j} A \Big) , 
\end{align*}
with $ {\bf w} = (w_1, w_2, \ldots, w_d) \in \RM^d$.
\end{thm}

\section{A vertex-face walk on a graph} 

We introduce a discrete-time quantum walk on an orientable embedding of a graph. 
An embedding of a graph is {\em circular} if every face is bounded by a cycle. 

Let $G$ be a graph, and ${\cal M} $ a circular embedding of $G$ on some orientable surface. 
Then we consider a consistent orientation of the faces, that is, for each edge $e$ shared by two faces $f$ and $h$, 
the direction $e$ receives in $f$ is opposite to the direction it receives in $h$. 
For such an orientation of faces, every arc belongs to exactly one face. 

Let $G$ have $n$ vertices, $m$ edges and $k$ faces, and $F(G)$ the set of faces of ${\cal M} $.  
Then, let ${\bf M} =( M_{ef} )_{e \in D(G); f \in F(G)} $ be the associated arc-face incidence matrix of $G$ defined as follows: 
\[
M_{ef} = \left\{
\begin{array}{ll}
1 & \mbox{if $e \in f$, } \\
0 & \mbox{otherwise. }
\end{array}
\right.
\] 
Furthermore, let ${\bf N} =( N_{ev} )_{e \in D(G); v \in V(G)} $ be the associated arc-origin incident matrix defined as follows: 
\[
N_{ev} = \left\{
\begin{array}{ll}
1 & \mbox{if $o(e)=v$, } \\
0 & \mbox{otherwise. }
\end{array}
\right.
\] 
Let $\hat{{\bf M}}$ and $\hat{{\bf N}}$ be normalized versions of ${\bf M}$ and ${\bf N}$, respectively. 
Note that 
\[
\hat{{\bf M}}^T \hat{{\bf M}} = {\bf I}_k, \ \  \  \  \hat{{\bf N}}^T \hat{{\bf N}} = {\bf I}_n . 
\]
Then the following unitary matrix ${\bf U} $ is the transition matrix of a {\em vertex-face walk} for ${\cal M} $:  
\[
{\bf U} =(2 \hat{{\bf M}} \hat{{\bf M}}^T - {\bf I}_{2m} )(2 \hat{{\bf N}} \hat{{\bf N}}^T - {\bf I}_{2m } ) . 
\]

Thus, we obtain the following formula for the transition matrix of a vertex-face walk on a graph.

\begin{thm}
Let $G$ be a circular embedding graph on some orientable surface, and have $n$ vertices, $m$ edges and $k$ faces. 
Then, for the transition matrix ${\bf U} $ of a vertex-face walk on $G$, 
\[
\det ( {\bf I}_{2m} -u {\bf U} )=(1-u )^{2m-n-k} (1+u )^{n-f} \det ((1+u)^2 {\bf I}_k 
-4u \hat{{\bf M}}^T \hat{{\bf N}} \hat{{\bf N}}^T \hat{{\bf M}} ). 
\]
\end{thm}

{\em Proof }.  By the definition of the transition matrix ${\bf U} $, we have  
\[
\begin{array}{rcl}
\  &   & \det ( {\bf I}_{2m} -u {\bf U} )= \det ( {\bf I}_{2m} 
-u(2 \hat{{\bf M}} \hat{{\bf M}}^T - {\bf I}_{2m} )(2 \hat{{\bf N}} \hat{{\bf N}}^T - {\bf I}_{2m } )) \\ 
\  &   &                \\ 
\  & = & 
\det ((1-u) {\bf I}_{2m} +2u \hat{{\bf N}} \hat{{\bf N}}^T -2u  \hat{{\bf M}} \hat{{\bf M}}^T (2 \hat{{\bf N}} \hat{{\bf N}}^T - {\bf I}_{2m } )) \\ 
\  &   &                \\ 
\  & = & (1-u)^{2m} \det ( {\bf I}_{2m} + \frac{2u}{1-u} \hat{{\bf N}} \hat{{\bf N}}^T 
- \frac{2u}{1-u} \hat{{\bf M}} \hat{{\bf M}}^T (2 \hat{{\bf N}} \hat{{\bf N}}^T - {\bf I}_{2m } )) \\ 
\  &   &                \\ 
\  & = & (1-u)^{2m} \det ( {\bf I}_{2m} - \frac{2u}{1-u} \hat{{\bf M}} \hat{{\bf M}}^T (2 \hat{{\bf N}} \hat{{\bf N}}^T - {\bf I}_{2m } )
( {\bf I}_{2m} + \frac{2u}{1-u} \hat{{\bf N}} \hat{{\bf N}}^T )^{-1} ) \\ 
\  &   &                \\ 
\  & \times & \det ( {\bf I}_{2m} + \frac{2u}{1-u} \hat{{\bf N}} \hat{{\bf N}}^T ) . 
\end{array}
\]

If ${\bf A}$ and ${\bf B}$ are an $m \times n $ and $n \times m$ 
matrices, respectively, then we have 
\[
\det ( {\bf I}_{m} - {\bf A} {\bf B} )= 
\det ( {\bf I}_n - {\bf B} {\bf A} ) . 
\]
Thus, we have 
\[
\begin{array}{rcl}
\  &   & \det ( {\bf I}_{2m} -u {\bf U} ) \\ 
\  &   &                \\ 
\  & = & (1-u)^{2m} \det ( {\bf I}_{k} - \frac{2u}{1-u} \hat{{\bf M}}^T (2 \hat{{\bf N}} \hat{{\bf N}}^T - {\bf I}_{2m } )
( {\bf I}_{2m} + \frac{2u}{1-u} \hat{{\bf N}} \hat{{\bf N}}^T )^{-1} \hat{{\bf M}} ) \\ 
\  &   &                \\ 
\  & \times & \det ( {\bf I}_{2m} + \frac{2u}{1-u} \hat{{\bf N}} \hat{{\bf N}}^T ) . 
\end{array}
\]

But, 
\[
\begin{array}{rcl}
\det ( {\bf I}_{2m} + \frac{2u}{1-u} \hat{{\bf N}} \hat{{\bf N}}^T ) 
& = & \det ( {\bf I}_{n} + \frac{2u}{1-u} \hat{{\bf N}}^T \hat{{\bf N}}) \\
\  &   &                \\ 
\  & = & \det ( {\bf I}_{n} + \frac{2u}{1-u} {\bf I}_n ) \\
\  &   &                \\ 
\  & = & (1+ \frac{2u}{1-u} )^n = \frac{(1+u)^n }{(1-u)^n } . 
\end{array}
\] 
Furthermore, we have 
\[
\begin{array}{rcl}
\  &   & ( {\bf I}_{2m} + \frac{2u}{1-u} \hat{{\bf N}} \hat{{\bf N}}^T )^{-1} \\ 
\  &   &                \\ 
\  & = & {\bf I}_{2m} - \frac{2u}{1-u} \hat{{\bf N}} \hat{{\bf N}}^T +( \frac{2u}{1-u} )^2 \hat{{\bf N}} \hat{{\bf N}}^T \hat{{\bf N}} \hat{{\bf N}}^T  
- ( \frac{2u}{1-u} )^3 \hat{{\bf N}} \hat{{\bf N}}^T \hat{{\bf N}} \hat{{\bf N}}^T \hat{{\bf N}} \hat{{\bf N}}^T+ \cdots \\
\  &   &                \\ 
\  & = & {\bf I}_{2m} - \frac{2u}{1-u} \hat{{\bf N}} \hat{{\bf N}}^T +( \frac{2u}{1-u} )^2 \hat{{\bf N}} \hat{{\bf N}}^T 
- ( \frac{2u}{1-u} )^3 \hat{{\bf N}} \hat{{\bf N}}^T + \cdots \\
\  &   &                \\ 
\  & = & {\bf I}_{2m} - \frac{2u}{1-u} (1 - \frac{2u}{1-u} +( \frac{2u}{1-u} )^2 - \cdots ) \hat{{\bf N}} \hat{{\bf N}}^T \\
\  &   &                \\ 
\  & = & {\bf I}_{2m} - \frac{2u}{1-u} /(1+ \frac{2u}{1-u} ) \hat{{\bf N}} \hat{{\bf N}}^T 
={\bf I}_{2m} - \frac{2u}{1+u} \hat{{\bf N}} \hat{{\bf N}}^T .
\end{array}
\]

Therefore, it follows that 
\[
\begin{array}{rcl}
\  &   & \det ( {\bf I}_{2m} -u {\bf U} ) \\ 
\  &   &                \\ 
\  & = & (1-u )^{2m} 
\det ( {\bf I}_{k} - \frac{2u}{1-u} \hat{{\bf M}}^T (2 \hat{{\bf N}} \hat{{\bf N}}^T - {\bf I}_{2m } ) 
( {\bf I}_{2m} - \frac{2u}{1+u} \hat{{\bf N}} \hat{{\bf N}}^T ) \hat{{\bf M}} )  \frac{(1+u)^n }{(1-u)^n } \\ 
\  &   &                \\ 
\  & = & (1-u )^{2m-n} (1+u)^n 
\det ( {\bf I}_{k} - \frac{2u}{1-u}  \hat{{\bf M}}^T (- {\bf I}_{2m} + \frac{2}{1+u} \hat{{\bf N}} \hat{{\bf N}}^T ) \hat{{\bf M}} ) \\ 
\  &   &                \\ 
\  & = & (1-u )^{2m-n} (1+u)^n 
\det ( {\bf I}_{k} + \frac{2u}{1-u} \hat{{\bf M}}^T \hat{{\bf M}} - \frac{4u}{1- u^2} \hat{{\bf M}}^T \hat{{\bf N}} \hat{{\bf N}}^T \hat{{\bf M}} ) \\ 
\  &   &                \\ 
\  & = & (1-u )^{2m-n-k} (1+u)^{n-k} 
\det ((1- u^2 ) {\bf I}_{k} +2u(1+u) {\bf I}_k -4u \hat{{\bf M}}^T \hat{{\bf N}} \hat{{\bf N}}^T \hat{{\bf M}} ) \\ 
\  &   &                \\ 
\  & = & (1-u )^{2m-n-k} (1+u)^{n-k} 
\det ((1+u )^2 {\bf I}_{k} -4u \hat{{\bf M}}^T \hat{{\bf N}} \hat{{\bf N}}^T \hat{{\bf M}} ) . 
\end{array}
\]
$\Box$

Substituting $u=1/ \lambda $, we obtain the following result.

\begin{cor}
Let $G$ be a circular embedding graph on some orientable surface, and have $n$ vertices, $m$ edges and $k$ faces. 
Then, for the transition matrix ${\bf U} $ of a vertex-face walk on $G$, 
\[
\det ( \lambda {\bf I}_{2m} - {\bf U} )=( \lambda -1)^{2m-n-k} ( \lambda +1 )^{n-f} \det (( \lambda +1)^2 {\bf I}_k 
-4 \lambda \hat{{\bf M}}^T \hat{{\bf N}} \hat{{\bf N}}^T \hat{{\bf M}} ) . 
\]
\end{cor}

{\em Proof }.  Let $u=1/ \lambda $. 
Then, by Theorem 4.1, we have 
\[
\det ({\bf I}_{2m} -1/  \lambda   {\bf U} )=(1-1/ \lambda )^{2m-n-k} (1+1/ \lambda )^{n-k} \det ((1+1/ \lambda )^2 {\bf I}_k -4/\lambda 
\hat{{\bf M}}^T \hat{{\bf N}} \hat{{\bf N}}^T \hat{{\bf M}} ) ,  
\]
and so, 
\[
\det ( \lambda  {\bf I}_{2m} - {\bf U} )=( \lambda -1)^{2m-n-k} ( \lambda +1 )^{n-f} \det (( \lambda +1)^2 {\bf I}_k 
-4 \lambda \hat{{\bf M}}^T \hat{{\bf N}} \hat{{\bf N}}^T \hat{{\bf M}} ) . 
\]
$\Box$

By Corollary 4.2, the following result holds. 
Let ${\rm Spec} ( {\bf F} )$ be the set of eigenvalues of a square matrix ${\bf F} $.

\begin{cor}
Let $G$ be a circular embedding graph on some orientable surface, and have $n$ vertices, $m$ edges and $k$ faces. 
Then the spectra of the transition matrix ${\bf U} $ are given as follows: 
\begin{enumerate} 
\item $2k$ eigenvalues: 
\[
\lambda =(2 \mu -1) \pm 2 \sqrt{ \mu ( \mu -1)} , \ \  \  \  \mu \in {\rm Spec} ( \hat{{\bf M}}^T \hat{{\bf N}} \hat{{\bf N}}^T \hat{{\bf M}} ) ; 
\]
\item $2m-n-k$ eigenvalues: 1; 
\item $n-k$ eigenvalues: -1. 
\end{enumerate} 
\end{cor}

{\em Proof }.  Let ${\rm Spec} ( \hat{{\bf M}}^T \hat{{\bf N}} \hat{{\bf N}}^T \hat{{\bf M}} )= \{ \lambda {}_1 , \ldots , \lambda {}_m \} $. 
Since $ \hat{{\bf M}}^T \hat{{\bf N}} \hat{{\bf N}}^T \hat{{\bf M}} $ is symmetric, we have 
\[
\lambda {}_1 , \ldots , \lambda {}_k \in \mathbb{R} . 
\]
Furthermore, by Corollary 4.2, we have 
\[
\begin{array}{rcl}
\  &   & \det ( \lambda  {\bf I}_{2m} - {\bf U} ) \\
\  &   &                \\ 
\  & = & ( \lambda -1 )^{2m-n-k} ( \lambda +1 )^{n-k} 
\prod_{ \mu \in {\rm Spec} ( \hat{{\bf M}}^T \hat{{\bf N}} \hat{{\bf N}}^T \hat{{\bf M}} )} (( \lambda +1)^2 -4 \mu \lambda ) \\  
\  &   &                \\ 
\  & = & ( \lambda -1 )^{2m-n-k} ( \lambda +1 )^{n-k} 
\prod_{ \mu \in {\rm Spec} ( \hat{{\bf M}}^T \hat{{\bf N}} \hat{{\bf N}}^T \hat{{\bf M}} )} ( \lambda {}^2 -2(2 \mu -1) \lambda +1) . 
\end{array}
\]

Solving $ \lambda {}^2 -2(2 \mu -1) \lambda +1=0$, we obtain 
\[
\lambda =(2 \mu -1) \pm 2 \sqrt{ \mu ( \mu -1)} . 
\] 
The result follows. 
$\Box$

\section{The vertex-face walk on the 2-dimensional finite torus}

At first, we state a result for the structure of the matrix $ \hat{{\bf M}}^T \hat{{\bf N}} \hat{{\bf N}}^T \hat{{\bf M}} $ by Zhan \cite{Zhan}. 

Let $G$ be a circular embedding graph on some orientable surface, and have $n$ vertices, $m$ edges and $k$ faces. 
For a face $f \in F(G)$, let $|f|$ be the number of vertices (or arcs) contained in the boundary of $f$.

\begin{pro}[Zhan]  
Let $G$ be a circular embedding graph on some orientable surface, and have $n$ vertices, $m$ edges and $k$ faces. 
Furthermore, let ${\bf K} = \hat{{\bf M}}^T \hat{{\bf N}} \hat{{\bf N}}^T \hat{{\bf M}} $. 
Then, for $f, h \in F(G)$, the $(f,h)$-entry of ${\bf K} $ is 
\[
K_{fh} = \frac{1}{ \sqrt{|f||h|}} \sum_{u \in f \cap h} \frac{1}{ \deg u} , 
\] 
where $f \cap h$ denotes the set of vertices used by both $f$ and $h$. 
\end{pro}

Now, for a natural number $N \geq 2$, let $G= T^2_N $ be the 2-dimensional finite torus (graph). 
Then we have 
\[
n=|V( T^2_N )|= N^2 , \ m=|E( T^2_N )|=2N^2 , \ k=|F( T^2_N )|= N^2 . 
\]
The 2-dimensional finite torus $T^2_N$ is circular embedded on the 2-dimensional torus $T^2 $, and its dual graph $(T^2_N )^*$ on $T^2 $ is isomorphic to 
itself. 
Furthermore, $T^2_N$ is a 4-regular graph, and the boundary of each face of $T^2_N$ has four vertices.

Then the following result holds for the transition matrix ${\bf U}$ of the vertex-face walk on $T^2_N $.

\begin{thm}  
Let $G= T^2_N $ be the 2-dimensional finite torus. 
Then 
\[
\det ( {\bf I}_{4N^2 } -u {\bf U} )=(1-u )^{2 N^2} \det ((1+u )^2 {\bf I}_{N^2 } - \frac{u}{4} ( {\bf A}^2 +2 {\bf A} )) . 
\]
\end{thm}

{\em Proof }.  Let $n=|V( T^2_N )|= N^2 , \ m=|E( T^2_N )|=2N^2 , \ k=|F( T^2_N )|= N^2 $, $V( T^2_N )= \{ v_1 , \ldots , v_n \} $, and 
$f_1 , \ldots , f_k \ (k=N^2 )$ be the faces of $T^2_N $. 
For $i=1, \ldots , k$, let $ {\bf f}_i =( (f_i)_u )_{u \in V(G)}$ be the $n$-dimensional vector such that 
\[
(f_i)_u = \left\{
\begin{array}{ll}
1 & \mbox{if $u \in f_i $, } \\
0 & \mbox{otherwise. }
\end{array}
\right.
\] 
Furthermore, for each $f \in F(G)$ and $u \in V(G)$, we have 
\[
|f|=4 \ \ and \ \ \deg u=4 . 
\]
By Proposition 5.1, we have 
\[
{\bf K} = \hat{{\bf M}}^T \hat{{\bf N}} \hat{{\bf N}}^T \hat{{\bf M}} = \frac{1}{16} ( {\bf f}_i \cdot {\bf f}_j )_{1 \leq i,j \leq k} ,  
\]
where ${\bf f}_i \cdot {\bf f}_j $ is the inner product of two vectors ${\bf f}_i $ and ${\bf f}_j $.  

But, we have 
\[
{\bf f}_i \cdot {\bf f}_j = |f_i \cap f_j | \ \  (1 \leq i,j \leq k) . 
\]
If $i=j$, then 
\[
{\bf f}_i \cdot {\bf f}_j =4 . 
\]
If $ |f_i \cap f_j |=2$, then $f_i $ and $ f_j $ are adjacent in the dual graph $G$. 
Since $T^2_N $ and $(T^2_N )^*$ are isomorphic, we have 
\[
{\bf A} ( G^* )= {\bf A} (G) . 
\]
Identifying $v_i $ and $f_i $ for each $i=1, \ldots , n$, we have 
\[
{\bf A} (G)_{ij} =1 \ \  \  if \ \ \  |f_i \cap f_j |=2 . 
\]

Next, if $|f_i \cap f_j |=1$, then $f_i $ and $ f_j $ are not adjacent in the dual graph $(T^2_N )^*$, and there exists 
a unique reduced path of length 2 from $f_i $ to $ f_j $ in $(T^2_N )^*$. 
Let ${\bf A}_2 (G)$ be the an $n \times n$ matrix such that 
\[
({\bf A}_2 (G))_{uv} = \left\{
\begin{array}{ll}
1 & \mbox{if there exists a reduced path from $u$ to $v$ in $G$, } \\
0 & \mbox{otherwise. }
\end{array}
\right.
\] 
Thus,  
\[
({\bf A}_2 (G))_{ij} =1 \ \ \ if \ \ \ |f_i \cap f_j |=1 . 
\]
Therefore it follows that 
\[
{\bf K} = \frac{1}{16} ( 4 {\bf I}_n +2 {\bf A} (G) + {\bf A}_2 (G)) . 
\]

But, since 
\[
{\bf A} (G)^2 = {\bf A}_2 (G) +4 {\bf I}_n ,  
\]
we have 
\[
{\bf K} = \frac{1}{16} ( 4 {\bf I}_n +2 {\bf A} (G)+ {\bf A}_2 (G)) = \frac{1}{16} ( 4 {\bf I}_n +2 {\bf A} (G) + {\bf A}^2 (G)-4 {\bf I}_n ) 
= \frac{1}{16} {\bf A}^2 + \frac{1}{8} {\bf A} .
\]
Thus, 
\[ 
\begin{array}{rcl}
\  &   & \det ( {\bf I}_{4N^2 } -u {\bf U} ) \\ 
\  &   &                \\ 
\  & = & (1-u )^{2m-n-k} (1+u)^{n-k} \det ((1+u )^2 {\bf I}_{k} -4u(  \frac{1}{16} {\bf A}^2 + \frac{1}{8} {\bf A} )) \\ 
\  &   &                \\ 
\  & = & (1-u )^{2 N^2} \det ((1+u )^2 {\bf I}_{N^2 } - \frac{u}{4} ( {\bf A}^2 +2 {\bf A} )) . 
\end{array}
\]
$\Box$

Therefore,

\begin{cor} 
Let $G= T^2_N $ be the 2-dimensional finite torus. 
Then  
\[ 
\begin{array}{rcl}
\  &   & \det ( \lambda {\bf I}_{4 N^2} - {\bf U} ) \\ 
\  &   &                \\ 
\  & = & ( \lambda -1 )^{2N^2} \prod^{N-1}_{k_1 =0} \prod^{N-1}_{k_2 =0} \left[ \lambda {}^2 - \lambda \{ \frac{1}{4} ( \cos \frac{2 \pi k_1 }{N} + \cos \frac{2 \pi k_2 }{N} )^2 
+ \cos \frac{2 \pi k_1 }{N} + \cos \frac{2 \pi k_2 }{N} -2 \} +1 \right] . 
\end{array}
\]
\end{cor}

{\em Proof }.  It is known that 
\[
{\rm Spec} ( T^2_N )= \{ 2 \cos \left( \frac{ 2 \pi k_1 }{N} \right) +2 \cos \left( \frac{2 \pi k_2 }{N} \right) \mid k_1 , k_2 =0,1, \ldots , N-1 \} . 
\]
By Theorem 5.2,  we have 
\[ 
\begin{array}{rcl}
\  &   & \det ( \lambda {\bf I}_{4 N^2} - {\bf U} ) \\ 
\  &   &                \\ 
\  & = & ( \lambda -1 )^{2N^2} \det (( \lambda +1 )^2 {\bf I}_{N^2 } - \frac{ \lambda }{4} ( {\bf A}^2 +2 {\bf A} )) \\
\  &   &                \\ 
\  & = & ( \lambda -1 )^{2N^2} \prod^{N-1}_{k_1 =0} \prod^{N-1}_{k_2 =0} (( \lambda +1)^2 - \frac{ \lambda }{4} (( \cos \frac{2 \pi k_1 }{N} + \cos \frac{2 \pi k_2 }{N} )^2 
+2 \cdot 2 ( \cos \frac{2 \pi k_1 }{N} + \cos \frac{2 \pi k_2 }{N} ))  \\
\  &   &                \\ 
\  & = & ( \lambda -1 )^{2N^2} \prod^{N-1}_{k_1 =0} \prod^{N-1}_{k_2 =0} ( \lambda {}^2 - \lambda ( \frac{1}{4} ( \cos \frac{2 \pi k_1 }{N} + \cos \frac{2 \pi k_2 }{N} )^2 
+ \cos \frac{2 \pi k_1 }{N} + \cos \frac{2 \pi k_2 }{N} -2)+1) . 
\end{array}
\]
$\Box$

As an application, we present the Walk/Zeta correspondence (see \cite{KomatsuEtAl2021b}) for the vertex-face walk on the 2-dimensional finite torus $G= T^2_N $. 

Similarly to the Walk/Zeta correspondence of \cite{KomatsuEtAl2021b}, we introduce a zeta function for the transition matrix ${\bf U} $ of the vertex-face walk on $G= T^2_N $ 
as follows: 
\[
\overline{\zeta } ( {\bf U} , T^2_N , u)= \det ( {\bf I}_{4 N^2} -u {\bf U} )^{-1/ N^2 } . 
\]
Let ${\bf A} = {\bf A} ( T^2_N )$. 
By Theorem 5.2 and Corollary 5.3, we have 
\[ 
\begin{array}{rcl}
\  &   & \overline{\zeta } ( {\bf U} , T^2_N , u)^{-1} = \det ( {\bf I}_{4 N^2} -u {\bf U} )^{1/ N^2 } \\ 
\  &   &                \\ 
\  & = & \{ (1-u )^{2 N^2} \det ((1+u )^2 {\bf I}_{N^2 } - \frac{u}{4} ( {\bf A}^2 +2 {\bf A} ))  \}^{1/ N^2 } \\
\  &   &                \\ 
\  & = & (1-u )^{2} \exp \left[ \frac{1}{N^2} \log \det \{ (1+u )^2 {\bf I}_{N^2 } - \frac{u}{4} ( {\bf A}^2 +2 {\bf A} ) \} \right] \\ 
\  &   &                \\ 
\  & = & (1-u )^{2} \\ 
\  &   &                \\ 
\  & \times & 
\exp  \left[ \frac{1}{N^2} \prod^{N-1}_{k_1 =0} \prod^{N-1}_{k_2 =0} \log \{ u^2 -u( \frac{1}{4} ( \cos \frac{2 \pi k_1 }{N} + \cos \frac{2 \pi k_2 }{N} )^2 
+ \cos \frac{2 \pi k_1 }{N} + \cos \frac{2 \pi k_2 }{N} -2)+1 \} \right] . 
\end{array}
\]
When $N \rightarrow \infty$, we have 
\[ 
\begin{array}{rcl}
\  &   & \lim_{N \rightarrow \infty} \overline{\zeta } ( {\bf U} , T^2_N , u)^{-1}\\ 
\  &   &                \\ 
\  & = & (1-u )^{2} 
\exp  \left[ \int^{2 \pi }_{0} \int^{2 \pi }_{0} \log \{ u^2 -u( \frac{1}{4} ( \cos \theta {}_1 + \cos \theta {}_2 )^2 
+ \cos \theta {}_1 + \cos \theta {}_2 -2)+1 \} \frac{d \theta {}_1 }{ 2 \pi } \frac{d \theta {}_2 }{ 2 \pi } \right] . 
\end{array}
\]

\begin{thm} 
Let $G= T^2_N $ be the 2-dimensional finite torus. 
Then 
\[ 
\begin{array}{rcl}
\  &   & \overline{\zeta } ( {\bf U} , T^2_N , u)^{-1} = (1-u )^{2} \\
\  &   &                \\ 
\  & \times & \exp \left[ \frac{1}{N^2} \prod^{N-1}_{k_1 =0} \prod^{N-1}_{k_2 =0} \log \{ u^2 -u( \frac{1}{4} ( \cos \frac{2 \pi k_1 }{N} 
+ \cos \frac{2 \pi k_2 }{N} )^2 + \cos \frac{2 \pi k_1 }{N} + \cos \frac{2 \pi k_2 }{N} -2)+1 \} \right]  
\end{array}
\]
and 
\[ 
\begin{array}{rcl}
\  &   & \lim_{N \rightarrow \infty} \overline{\zeta } ( {\bf U} , T^2_N , u)^{-1} \\ 
\  &   &                \\ 
\  & = & (1-u )^{2} \exp \left[ \int^{2 \pi }_{0} \int^{2 \pi }_{0} \log \{ u^2 -u( \frac{1}{4} \{ \cos \theta {}_1 + \cos \theta {}_2 )^2 
+ \cos \theta {}_1 + \cos \theta {}_2 -2)+1 \} \frac{d \theta {}_1 }{ 2 \pi } \frac{d \theta {}_2 }{ 2 \pi } \right] . 
\end{array}
\]
\end{thm}

\begin{small}
\bibliographystyle{jplain}

\end{small}

\end{document}